\documentstyle[12pt,amsmath,amssymb]{article}
\begin{document}

\newtheorem{lem}{Lemma}[section]
\newtheorem{theorem}{Theorem}[section]
\newtheorem{prop}{Proposition}[section]
\newtheorem{rem}{Remark}[section]
\newtheorem{define}{Definition}[section]
\newtheorem{cor}{Corollary}[section]
\allowdisplaybreaks

\makeatletter\@addtoreset{equation}{section}\makeatother
\def\theequation{\arabic{section}.\arabic{equation}}

\newcommand{\D}{{\cal D}}

\newcommand{\N}{{\Bbb N}}
\newcommand{\C}{{\Bbb C}}
\newcommand{\Z}{{\Bbb Z}}
\newcommand{\R}{{\Bbb R}}
\newcommand{\la}{\langle}
\newcommand{\ra}{\rangle}
\newcommand{\rom}[1]{{\rm #1}}
\newcommand{\FC}{{\cal F}C_{\mathrm b}(C_0(X),\Gamma)}
\newcommand{\eps}{\varepsilon}
\newcommand{\dd}{\overset{{.}{.}}}
\newcommand{\fii}{\varphi}

\newcommand{\supp}{\operatorname{supp}}

\def\stackunder#1#2{\mathrel{\mathop{#2}\limits_{#1}}}
\newcommand{\FCo}{{\cal F}C_{\mathrm b}^\infty({\cal D},\dd\Gamma)}

\newcommand{\EG}{{\cal E}_{\mathrm G}}
\newcommand{\EK}{{\cal E}}

\renewcommand{\author}[1]{\medskip{\large #1}\par\medskip}
\begin{center}{\Large \bf
  Equilibrium  Kawasaki dynamics of continuous particle systems
}\end{center}

{\large Yuri Kondratiev}\\
 Fakult\"at f\"ur Mathematik, Universit\"at
Bielefeld, Postfach 10 01 31, D-33501 Bielefeld, Germany;
 BiBoS, Univ.\ Bielefeld,
Germany.\\ e-mail:
\texttt{kondrat@mathematik.uni-bielefeld.de}\vspace{2mm}

{\large Eugene Lytvynov}\\ Department of Mathematics,
University of Wales Swansea, Singleton Park, Swansea SA2 8PP, U.K.\\
e-mail: \texttt{e.lytvynov@swansea.ac.uk}\vspace{2mm}

{\large Michael R\"ockner}\\
 Fakult\"at f\"ur Mathematik, Universit\"at
Bielefeld, Postfach 10 01 31, D-33501 Bielefeld, Germany;
 BiBoS, Univ.\ Bielefeld,
Germany.\\ e-mail: \texttt{roeckner@mathematik.uni-bielefeld.de }

{\small

\begin{center}
{\bf Abstract}
\end{center}

\noindent  We construct a new  equilibrium dynamics of infinite
particle systems in a Riemannian manifold $X$. This dynamics is an
analog of the  Kawasaki dynamics of lattice spin systems. The
Kawasaki  dynamics now is a process where interacting particles
randomly hop over $X$.  We establish conditions on the {\it a priori\/}
explicitly
 given symmetrizing  measure and
 the generator of this dynamics, under which
 a corresponding conservative Markov processes exists. We also
 outline two types of scaling limit of the equilibrium Kawasaki dynamics: one
 leading to an equilibrium Glauber dynamics in continuum (a birth-and-death
 process), and the other leading to a diffusion dynamics of
 interacting particles (in particular, the gradient stochastic
 dynamics).

 } \vspace{2mm}

 \vspace{1.5mm}

  \vspace{1.5mm}

\section{Introduction}

In the classical $d$-dimensional Ising model with   spin space
$S=\{-1,1\}$, the Kawasaki  dynamics means that  pairs of
neighboring particles with different spins randomly exchange their
spin values. The generator of this dynamics is given by
$$ (Hf)(\sigma)=\sum_{x\in\Z^d}\,\sum_{y\in\Z^d,\, |x-y|=1}c(x,y,\sigma)(\nabla_{xy}f)(\sigma),$$
where $$ (\nabla_{xy}f)(\sigma)=f(\sigma^{xy})-f(\sigma),$$
$\sigma^{xy}$ denoting the configuration  $\sigma$ in which the
particles at sites $x$ and $y$ have exchanged their spin values.
Under appropriate conditions on the coefficient  
$c(x,y,\sigma)$, the  corresponding dynamics has a Gibbs measure as
symmetrizing (hence invariant) measure. We refer, e.g., to
\cite{Liggett} for a discussion of the  Kawasaki dynamics of
lattice spin systems.

Let us now interpret a lattice system with spin space $S=\{-1,1\}$
as a model of a lattice gas. Then $\sigma(x)=1$ means that there is
a particle at site $x$, while $\sigma(x)=-1$  means that the site
$x$ is empty. The Kawasaki dynamics of such a system means that
particles randomly hop from one site to another.

If we consider a continuous particle system, i.e., a system of
particles which can take any position in the Euclidean  space
$\R^d$, then an analog of  the Kawasaki dynamics should be a process
in which   particles randomly hop over the space $\R^d$.
 The generator of such a process is  informally given by
 \begin{equation}\label{hjcdf}
(HF)(\gamma)=\sum_{x\in\gamma}\int_{\R^d}
c(x,y,\gamma)(D^{-+}_{xy}F)(\gamma)\,dy,\end{equation} where
\begin{equation}\label{222}
(D_{xy}^{-+}F)(\gamma)=F(\gamma\setminus x\cup
y)-F(\gamma)\end{equation} and the coefficient $c(x,y,\gamma)$
describes the rate at which the particle $x$ of the configuration
$\gamma$ jumps to $y$. Here and below, for simplicity of notations,
we just write $x$ instead of $\{x\}$.

In \cite{G2}, Gl\"otzl  considered the  formal generator
\eqref{hjcdf} and derived conditions on the coefficient
$c(x,y,\gamma)$ under which the operator \eqref{hjcdf} becomes 
symmetric  in the space $L^2(\mu)$, where $\mu$ is a given
Gibbs measure. However,  till now, there has  been no proof of the
very existence of a Kawasaki dynamics of an infinite system of
interacting particles  in continuum.

Thus, the aim of this paper is to present a general theorem  on the
existence of an equilibrium Kawasaki dynamics of  a continuous
particle system, which  has  a Gibbs measure as  symmetrizing (hence
invariant) measure. We shall also  consider some examples of such a
dynamics. Finally, we shall
 outline two types of scaling limit of the equilibrium Kawasaki dynamics: one
 leading to an equilibrium Glauber dynamics in continuum (a birth-and-death
 process), and the other leading to a diffusion dynamics of
 interacting particles (in particular, the gradient stochastic
 dynamics). A detailed study of these scaling limits will be given elsewhere
 \cite{FKL,KKL}.

Let us now briefly outline the structure of the paper.

In Section~\ref{fsc}, we fix a Riemannian manifold $X$  as  underlying space  (the position space of the particles)
and the space $\Gamma$ of all locally finite configurations in $X$.
The restriction  to the Riemannian manifold case is mainly motivated by the necessity to have constructive conditions for the existence of equilibrium states for interacting particle systems in $X$. Let us stress that all general statements of the paper (with minor changes) remain valid for much more general underlying spaces.

We next recall the definition of a Gibbs measure $\mu$  on $\Gamma$
which corresponds to a relative energy  $E(x,\gamma)$ of the
interaction between a particle $x$ and a configuration $\gamma$.
About the measure $\mu$ we assume that it has correlation functions
which satisfy the classical Ruelle bound. We also  present some
examples of a Gibbs measure corresponding to  a pair potential
$\phi$. It should be mentioned   that, although in   the examples we
deal with a pair potential $\phi$, our general theory for existence
of  dynamics holds for a general relative energy $E(x,\gamma)$.

  Next, in Section~\ref{rfyje},   under mild conditions  on $E(x,\gamma)$,
  we prove that there exists a  Hunt processes  ${\bf M}$
  on  $\Gamma$ which is properly associated with the Dirichlet form
of the  Kawasaki dynamics. In particular, ${\bf M}$  is a
conservative Markov process on $\Gamma$ with {\it cadlag} paths, and
has $\mu$ as symmetrizing, hence invariant  measure. We also
characterize this process in terms of the corresponding martingale
problem. Furthermore, we discuss the explicit form of the
$L^2(\mu)$-generator of this process on the set of continuous
bounded cylinder functions. 
In this section, we use the theory of Dirichlet forms \cite{MR}, and
in particular, some ideas and techniques developed in
\cite{KL,KLR2004,MR98,RS98}.

In Section~\ref{lkhuhy} we consider some examples of Kawasaki
dynamics. 

Finally, in Section~\ref{newsection}, by analogy with the Kawasaki
dynamics, we formulate conditions which guarantee the existence of
an equilibrium  Glauber (birth-and-death) dynamics in continuum
(compare with \cite{HS,KL,P}). We then
 outline  the above mentioned  scaling limits of the equilibrium Kawasaki dynamics.

We conclude this section with the following remarks. In a bounded domain, a Kawasaki dynamics
can be described as a jump Markov process. However, in the infinite volume, this  dynamics does not belong to this class, since in any time interval $[0,t]$, the dynamics has an infinite number of jumps.
Note also that 
 the set of symmetrizing measures of a given  Kawasaki dynamics consists of all  grand-canonical 
 Gibbs measures corresponding to a given relative energy of interaction and   {\it any} activity parameter $z>0$. This fact makes it especially interesting to study the hydrodynamic behavior  of the Kawasaki dynamics, cf.\ \cite{777,323232}.  
Finally, note   a similarity between the Kawasaki dynamics and the diffusion dynamics of continuous particle systems, e.g.\
 \cite{111,121212,KLR2004,MR98}.  Namely, both types of dynanics have (at least heuristically) conserved particle numbers and the same set of symmetrizing measures. Therefore, just as in the diffusion case, it is natural to study the scaling limit of equilibrium fluctuations for the Kawasaki dynamics, which is the subject of \cite{171717}.   

\section{Gibbs measures on configuration spaces}\label{fsc}

Let $X$ be a connected oriented $C^\infty$ manifold.  We denote the Riemannian distance on $X$ by  $\operatorname{dist}$. Let ${\cal B}(X)$ denote the Borel $\sigma$-algebra on $X$ and   $m$  the volume measure on $X$.

The configuration space $\Gamma:=\Gamma_{X}$ over $X$ is defined as the set of all subsets of $X$ which are
locally finite: $$\Gamma:=\big\{\,\gamma\subset X : \,
|\gamma_\Lambda|<\infty\text{ for each compact }\Lambda\subset
X\,\big\},$$ where $|\cdot|$ denotes the cardinality of a set
and $\gamma_\Lambda:= \gamma\cap\Lambda$. One can identify any
$\gamma\in\Gamma$ with the positive Radon measure
$\sum_{x\in\gamma}\eps_x\in{\cal M}(X)$,  where  $\eps_x$ is
the Dirac measure with mass at $x$,
$\sum_{x\in\varnothing}\varepsilon_x{:=}$zero measure, and ${\cal
M}(X)$
 stands for the set of all
positive  Radon  measures on   ${\cal
B}(X)$. The space $\Gamma$ can be endowed with the relative
topology as a subset of the space ${\cal M}(X)$ with the vague
topology, i.e., the weakest topology on $\Gamma$ with respect to
which  all maps $$\Gamma\ni\gamma\mapsto\la
f,\gamma\ra:=\int_{X} f(x)\,\gamma(dx) =\sum_{x\in\gamma}f(x),\qquad
f\in C_0(X),$$ are continuous. Here, $C_0(X)$ is the space
of all continuous  real-valued functions on $X$ with compact
support. We shall  denote  the Borel
$\sigma$-algebra on $\Gamma$ by ${\cal B}(\Gamma)$.

Now we proceed to consider Gibbs measures on $\Gamma$. For $\gamma\in\Gamma$ and $x\in X$, we consider a relative energy $E(x,\gamma)$ of interaction between a particle located at $x$ and the configuration $\gamma$. We suppose that the mapping $E$ is measurable and $E(x,\gamma)\in(-\infty,+\infty]$.

A probability measure $\mu$ on $(\Gamma,{\cal B}(\Gamma))$ is called a (grand-canonical) Gibbs measure corresponding to activity $z>0$ and the relative energy $E$
if it satisfies the Georgii--Nguyen--Zessin identity (\cite[Theorem~2]{NZ}, see also
\cite[Theorem~2.2.4]{KunaPhD}):
\begin{equation}\int_\Gamma \mu(d\gamma)\int_{X} \gamma(dx)
F(x,\gamma) =\int_\Gamma \mu(d\gamma)\int_{X}
zm(dx)\exp\left[-E(x,\gamma)\right] F(x,\gamma\cup
x)\label{fdrtsdrt}\end{equation} for any measurable function
$F:X\times \Gamma \to[0,+\infty]$. Let ${\cal G}(z,E)$
denote the set of all Gibbs measures corresponding to $z$ and
$E$.

In particular, if $E(x,\gamma)\equiv0$, then \eqref{fdrtsdrt} is
the Mecke identity, which holds if and only if $\mu$ is the
Poisson measure $\pi_z$ with intensity measure $z m(dx)$.

We  assume that
 \begin{equation}\label{jkhgiu}
E(x,\gamma)\in\R\quad\text{for $m\otimes \mu $-a.e.\ $(x,\gamma)\in X\times\Gamma$}.\end{equation} Furthermore, we assume that,
for any $n\in\N$\rom, there exists a non-negative measurable
symmetric function $k_\mu^{(n)}$ on $X^n$
 such
that, for any measurable symmetric function
$f^{(n)}:X^n\to[0,\infty]$,
\begin{align} &\int_\Gamma \sum_{\{x_1,\dots,x_n\}\subset\gamma}
f^{(n)}(x_1,\dots,x_n)\,\mu(d\gamma)\notag\\&\qquad =\frac1{n!}\,
\int_{X^n} f^{(n)}(x_1,\dots,x_n)
k_\mu^{(n)}(x_1,\dots,x_n)\,m(dx_1)\dotsm
m(dx_n),\label{6t565r7}\end{align} and
\begin{equation}\label{swaswea957}\forall (x_1,\dots,x_n)\in X^n:\quad
k_\mu^{(n)}(x_1,\dots,x_n)\le \xi^n,\end{equation} where $\xi> 0$
is independent of $n$.
The functions $k_\mu^{(n)}$, $n\in\N$, are called the  correlation
functions of the measure $\mu$, while \eqref{swaswea957} is called
the Ruelle bound.

Notice that any probability  measure $\mu$ on $(\Gamma,{\cal B}(\Gamma))$
satisfyng the Ruelle bound  has all local moments finite, i.e.,
\begin{equation}\label{sweqaw}\int_\Gamma \la f,\gamma\ra^n\,\mu(d\gamma)<\infty,\qquad
f\in C_0(X),\ f\ge0,\  n\in\N.\end{equation}

Let us give examples of a Gibbs measure corresponding to a pair potential $\phi$ and satisfying  the above assumptions.

Let $\phi:X^2\to(-\infty,+\infty]$ be a symmetric measurable function such that $\phi(x,y)\in\R$ for any $x,y\in X$, $x\ne y$. For each $x\in X$ and $\gamma\in\Gamma$, we define $$
E(x,\gamma):=\begin{cases}
\sum_{y\in\gamma}\phi(x,y),&\text{if } \sum_{y\in\gamma}|\phi(x,y)|<\infty,\\ +\infty,&\text{otherwise}.
\end{cases}$$

Let us formulate some conditions on the pair potential $\phi$.

\begin{description}

\item[(S)] ({\it Stability})
 There exists $B\ge0$ such that, for any $\gamma\in\Gamma$, $|\gamma|<\infty$,
$$\sum_{\{x,y\}\subset\gamma}\phi(x,y)\ge
-B|\gamma|.$$

\item[(I)] ({\it Integrability})
We have $$C:=\sup_{x\in X}\int _X
|\exp[-\phi(x,y)]-1|\,m(dy)<\infty.$$

\item[(F)] ({\it Finite range}) There exists $R>0$ such that $$
\phi(x,y)=0\quad\text{if }\operatorname{dist}(x,y)\ge R.$$

\end{description}

Note that if $\phi$ satisfies (F), then $E(x,\gamma)\in\R$ for any $\gamma\in\Gamma$  and $x\in X\setminus
\gamma$.

\begin{theorem}[\rom{\cite{KKS98,KunaPhD,K}}]\label{gdtrsresa}

\rom{1)} Let \rom{(S), (I),} and \rom{(F)} hold\rom, and let $z>0$
be such
 that $$z<\frac1{2e}\,( e^{2 B}C)^{-1},$$ where
$B$ and $C$ are as in \rom{(S)} and \rom{(I),} respectively\rom.
Then there exists a Gibbs measure $\mu\in{\cal G}(z,E)$
whose correlation functions $k_\mu^{(n)}$  exist and satisfy the Ruelle bound\rom.

\rom{2)} Let $\phi$ be a non-negative potential which fulfills
\rom{(I)} and \rom{(F).} Then for each $z>0$\rom, there
exists a Gibbs measure $\mu\in{\cal G}(z,E)$ whose
correlation functions $k_\mu^{(n)}$  exist and satisfy
the Ruelle bound.

\end{theorem}

Assume  now that   $X=\R^d$, $d\in\N$,
and assume that $\phi$ is  translation invariant, i.e., $\phi(x,y)=\tilde \phi(x-y)$,
where $\tilde\phi:\R\to(-\infty,\infty]$ is such that $\tilde\phi(x)\in\R$ for $x\ne0$  and $\tilde\phi(-x)=\tilde\phi(x)$
for all $x\in\R^d$.  In this case, the conditions on $z$ and $\phi$ can be significantly weakened.
First, we note that the condition (I) now looks as follows:
$$C:=\int _{\R^d}
|\exp[-\tilde\phi(x)]-1|\,m(dx)<\infty.$$

For the notion of a superstable,  lower regular potential and  the notion of a tempered Gibbs measure, appearing in the following theorem, see \cite{Ru70}.

\begin{theorem}[\rom{\cite{Ru69,Ru70}}]
\label{yte6te645}

 Assume that $X=\R^d$ and $\phi$ is translation invariant.

\rom{1)} Let \rom{(S)} and \rom{(I)} hold and let $z>0$
be such
 that $$z<\frac1{e}\,( e^{2 B}C)^{-1},$$ where
$B$ and $C$ are as in \rom{(S)} and \rom{(I),} respectively\rom.
Then there exists a Gibbs measure $\mu\in{\cal G}(z,E)$
whose  correlation functions exist and satisfy the Ruelle bound\rom.

\rom{2)} Let $\phi$ be a non-negative potential which fulfills
\rom{(I).} Then\rom, for each $z>0$\rom, there
exists a Gibbs measure $\mu\in{\cal G}(z,E)$ whose correlation functions exist and satisfy
the Ruelle bound.

\rom{3)} Let $\phi$ satisfy \rom{(I)} and additionally let $\phi$ be a superstable\rom, lower regular potential.
Then the set ${\cal G}_{{\mathrm temp}}(z,E)$ of all tempered  Gibbs measures is non-empty and each measure from
${\cal G}_{{\mathrm temp}}(z,E)$ has correlation functions which satisfy the Ruelle bound\rom.

\end{theorem}

We also have the following lemma, which follows from (the proof of) \cite[Lemma 3.1]{KLR2004}.

\begin{lem}\label{lohkjj}
Let $X=\R^d$ and  let  $\phi$\rom, $z$\rom, and $\mu\in{\cal G}(z,E)$ be as in one of the statements of Theorem~\rom{\ref{yte6te645}.}
Assume additionally that there exists $r>0$ such that \begin{equation}
\label{fytr} \sup_{x\in B(r)^c}\tilde\phi(x)<\infty, \end{equation} where $B(r)$ denotes the ball in $\R^d$ of radius $r$ centered at the origin\rom. Then  \eqref{jkhgiu} holds\rom.

\end{lem}

\section{Existence  results}\label{rfyje}

In what follows, we shall consider a Gibbs measure $\mu\in{\cal G}(z,E)$ as in Section~\ref{fsc}, i.e., a probability measure $\mu$ on $(\Gamma,{\cal B}(\Gamma))$ which satisfies \eqref{fdrtsdrt}--\eqref{swaswea957}. We introduce the set $\FC$
of all functions of the form $$\Gamma\ni \gamma\mapsto F(\gamma)=g_F(\la\varphi_1,\gamma\ra,\dots , \la\varphi_N,\gamma\ra), $$ where $N\in\N$, $\varphi_1,\dots,\varphi_N\in C_0(X)$ and $g_F\in C_{\mathrm b}(\R^N)$,
where $C_{\mathrm b}(\R^N)$ denotes the set of all continuous bounded functions on $\R^N$.

We consider a measurable mapping $$ X\times X\times\Gamma\ni
(x,y,\gamma)\mapsto c(x,y,\gamma)\in[0,\infty).$$ We assume that,
 for each compact $\Lambda\subset X$,
\begin{equation}
\int_\Gamma\mu(d\gamma)\int_X\gamma(dx)\int_X
m(dy)c(x,y,\gamma)(\pmb1_\Lambda(x)+\pmb1_\Lambda(y))<\infty,
\label{gyu}\end{equation} where $\pmb1_\Lambda$ denotes the
indicator of $\Lambda$.

For each function $F:\Gamma\to\R$, $\gamma\in\Gamma$, and $x,y\in
X$, we recall the notation  \eqref{222}. Then we define a bilinear form
\begin{equation} {\cal E}
(F,G):= \int_\Gamma \mu(d\gamma)\int_X \gamma(dx)\int_X zm(dy)
c(x,y,\gamma) (D_{xy}^{-+}F)(\gamma)
(D_{xy}^{-+}G)(\gamma),\label{2}
\end{equation}
where $F,G\in\FC$. Below we shall  show that ${\cal E}$ corresponds
to a  Kawasaki dynamics.

We note that, for any $F\in\FC$, there exist a compact
$\Lambda\subset X$ and $C_1>0$ such that
$$
|(D_{xy}^{-+}F)(\gamma)|\le C_1
(\pmb1_\Lambda(x)+\pmb1_\Lambda(y)),\qquad \gamma\in\Gamma,\ x,y\in
X.$$ Therefore, by  \eqref{gyu}, the right hand sides of formula
 \eqref{2} is well-defined and finite.

\begin{lem} We have ${\cal E}(F,G)=0$ for all $F,G\in\FC$ such that
$F=0$ $\mu$-a\rom.e\rom.
\end{lem}

\noindent {\it Proof.}  It suffices to show that, for $F\in\FC$, $F=0$ $\mu$-a.e., we have
$(D_{x,y}^{-+}F)(\gamma)=0$ $\tilde\mu$-a.e., where $\tilde\mu$ is the measure on $X\times X\times\Gamma$ defined by \begin{equation}\label{gyugy}\tilde
\mu(dx,dy,d\gamma):=\gamma(dx)\, zm(dy)\,\mu(d\gamma).\end{equation}

Let $\Lambda$ be a compact subset of $X$.  We have:
$$ \int_\Gamma \mu(d\gamma)\int_\Lambda \gamma(dx)\int_\Lambda zm(dy)|F(\gamma)|=
 \int_\Gamma \mu(d\gamma) |F(\gamma)|\int_\Lambda \gamma(dx)\int_\Lambda zm(dy)=0, $$
which implies that $F(\gamma)=0$ $\tilde\mu$-a.e. Next, by \eqref{fdrtsdrt} and \eqref{jkhgiu},
\begin{multline}\label{yuuy} \int_\Gamma \mu(d\gamma)\int_\Lambda\gamma(dx)\int_\Lambda zm(dy)|F(\gamma\setminus x\cup y)|\\
=
\int_\Gamma\mu(d\gamma)|F(\gamma)| \int_\Lambda\gamma(dx)\int_\Lambda zm(dy)
\exp[-E(y,\gamma)+
E(x,\gamma\setminus x\cup y)].\end{multline}
Since $F$ is bounded, by \eqref{sweqaw}, the integrals in \eqref{yuuy} are finite. Therefore, \begin{equation}\label{ghiuygiyg}|F(\gamma)|\exp[-E(y,\gamma)+
E(x,\gamma\setminus x\cup y)]<\infty \quad \text{for $\tilde \mu$-a.e.\ $(x,y,\gamma)\in X\times X\times\Gamma$}.\end{equation} Since $F=0$ $\mu$-a.e., by \eqref{yuuy} and \eqref{ghiuygiyg}, $F(\gamma\setminus x\cup y)=0$ $\tilde\mu$-a.e.\quad $\square$\vspace{2mm}

Thus, $({\cal E},\FC)$ is a well-defined bilinear form on
$L^2(\Gamma,\mu)$.

\begin{lem}\label{uiwfg}
The bilinear form $({\cal E},\FC)$ is closable on $L^2(\Gamma,\mu)$
and its closure will be denoted by $({\cal E},D({\cal E}))$\rom.
\end{lem}

\noindent {\it Proof.} Let $(F_n)_{n=1}^\infty$ be a sequence in $\FC$ such that
$\|F_n\|_{L^2(\mu)}\to0$ as $n\to\infty$ and
\begin{equation}\label{rdd}\EK(F_n-F_k)\to0\quad \text{as $n,k\to\infty$}.\end{equation} Here and below, $\EK(F)$ stays for $\EK(F,F)$.
To prove the closability of $\EK$ it suffices to show that there exists a subsequence $(F_{n_k})_{k=1}^\infty$ such that
$\EK(F_{n_k})\to0$ as $k\to\infty$.

Let $\Lambda$ be a compact subset of $X$. By \eqref{sweqaw}, we have $$ \int_\Gamma\mu(d\gamma)\int_\Lambda\gamma(dx)\, |F_n(\gamma)|\le \|F_n\|_{L^2(\mu)}\left(\int_\Gamma \la \pmb1_\Lambda,\gamma\ra^2\,\mu(d\gamma)\right)^{1/2}\to0\quad\text{as }n\to\infty.$$ Therefore, there exists a subsequence of $(F_n)_{n=1}^\infty$, denoted by $(F_n^{(1)})_{n=1}^\infty$, such that $F_n^{(1)}(\gamma)\to0$ for $\gamma(dx)\mu(d\gamma)$-a.e.\ $(x,\gamma)\in \Lambda\times\Gamma$. Hence, there exists a subsequence $(F_n^{(2)})_{n=1}^\infty $ of  $(F_n^{(1)})_{n=1}^\infty$ such that
$F_n^{(2)}(\gamma)\to0$ for $\gamma(dx)\mu(d\gamma)$-a.e.\ $(x,\gamma)\in X\times\Gamma$.

Next, analogously to  \eqref{yuuy},
\begin{align*}&
\int_\Gamma \mu(d\gamma) \int_\Lambda \gamma(dx)\int_\Lambda zm(dy)\, \exp[-E(y,\gamma)+E(x,\gamma\setminus x\cup y)]\,|F_n^{(2)}(\gamma\setminus x\cup y)|\\&\qquad =\int_\Gamma \mu(d\gamma)\int_\Lambda zm(dx)\int_\Lambda
\gamma(dy)|F_n(\gamma)|\\&\qquad\le \|F_n^{(2)}\|_{L^2(\mu)} zm(\Lambda) \left(\int_\Gamma \la \pmb1_\Lambda,\gamma\ra^2\,\mu(d\gamma)\right)^{1/2}\to0\quad\text{as }n\to\infty.
\end{align*}
By virtue of \eqref{jkhgiu}, $$\exp[-E(y,\gamma)+E(x,\gamma\setminus x\cup y)]\in(0,+\infty]\quad \text{for $\tilde \mu$-a.e.\
$(x,y,\gamma)\in X\times X\times \Gamma$}.$$
Therefore, there exists a subsequence $(F_n^{(3)})_{n=1}^\infty$ of
$(F_n^{(2)})_{n=1}^\infty $ such that $F_n^{(3)}(\gamma\setminus x\cup y)\to0$ for $\tilde\mu$-a.e.\
$(x,y,\gamma)\in X\times X\times\Gamma$, where the measure $\tilde \mu$ is defined by \eqref{gyugy}.

Thus, \begin{equation}\label{hyugfu} (D_{xy}^{-+}F_n^{(3)})(\gamma)\to0\quad\text{as $n\to\infty$ for $\tilde\mu$-a.e.\
$(x,y,\gamma)\in X\times X\times\Gamma$}.\end{equation}
Now, by \eqref{hyugfu} and Fatou's lemma
\begin{align*}
\EK(F_n^{(3)})&=\int c(x,y,\gamma)(D_{xy}^{-+}F_n^{(3)})(\gamma)^2\,\tilde\mu(dx,dy,d\gamma)\\&=
\int c(x,y,\gamma)\left((D_{xy}^{-+}F_n^{(3)})(\gamma)-\lim_{m\to\infty}(D_{xy}^{-+}F_m^{(3)})(\gamma)\right)^2\,\tilde\mu(dx,dy,d\gamma)\\&\le
\liminf_{m\to\infty}\int c(x,y,\gamma)((D_{xy}^{-+}F^{(3)}_n)(\gamma)-(D_{xy}^{-+}F^{(3)}_m)(\gamma))^2\,\tilde\mu(dx,dy,d\gamma)\\&=\liminf_{m\to\infty}\EK(F_n^{(3)}-F_m^{(3)}),
\end{align*}
which by \eqref{rdd} can be made arbitrarily small for $n$ large enough.\quad $\square$\vspace{2mm}

For the notion of a Dirichlet form,  appearing in the
following lemma, we refer to e.g.\ \cite[Chap.~I, Sect.~4]{MR}.

\begin{lem}\label{guzazagus}
$({\cal E},D({\cal E})$ is a  Dirichlet form on $L^2(\Gamma,\mu)$\rom.
\end{lem}

\noindent {\it Proof.} On $D(\EK)$ we consider the norm $\|F\|_{D(\EK)}:=(\|F\|^2_{L^2(\mu)}+\EK(F))^{1/2}$, $F\in D(\EK)$.
For any $F,G\in\FC$, we define $$ S(F,G)(x,y,\gamma):=c(x,y,\gamma)(D_{xy}^{-+}F)(\gamma) (D_{xy}^{-+}G)(\gamma)
,\qquad x,y\in X,\gamma\in\Gamma.$$
Using the Cauchy  inequality, we conclude that $S$ extends to a bilinear continuous map from $(D(\EK),
\|\cdot\|_{D(\EK)})\times (D(\EK),
\|\cdot\|_{D(\EK)})$ into $L^1(X\times X\times\Gamma,\tilde\mu)$. Let $F\in D(\EK)$ and consider any sequence $(F_n)_{n=1}^\infty$ in $\FC$ such that $F_n\to F$ in $(D(\EK),\|\cdot\|_{D(\EK)})$.
In particular, $F_n\to F$ in $L^2(\mu)$. Then, analogously to the proof of Lemma~\ref{uiwfg}, for some subsequence
$(F_{n_k})_{k=1}^\infty $, we get $$
(D_{xy}^{-+}F_{n_k})(\gamma)\to (D_{xy}^{-+}F)(\gamma)\quad \text{for $\tilde\mu$-a.e.\ }(x,y,\gamma)\in X\times X\times\Gamma. $$ Therefore, for any $F,G\in D(\EK)$,
\begin{equation}\label{iuedw} S(F,G)(x,y,\gamma):=c(x,y,\gamma)(D_{xy}^{-+}F)(\gamma) (D_{xy}^{-+}G)(\gamma)\quad
\text{for $\tilde\mu$-a.e.\ }(x,y,\gamma)\in X\times X\times\Gamma \end{equation}
and \begin{equation}\label{guyiuhiu}\EK(F,G)=\int S(F,G)(x,y,\gamma)\,\tilde\mu(dx,dy,d\gamma).\end{equation}

Define $\R\ni x\mapsto g(x){:=}(0\vee x)\wedge 1$.  We again fix
any $F\in D(\EK)$ and  let $(F_n)_{n=1}^\infty$ be a sequence
of functions from $\FC$ such that $F_n\to F$ in $(D(\EK),\|\cdot\|_{D(\EK)})$. Consider the sequence
$(g(F_n))_{n\in\N}$. We evidently have: $g(F_n)\in\FC$ for each
$n\in\N$ and, by the dominated   convergence theorem,
 $g(F_n)\to g(F)$ as $n\to\infty$ in
$L^2(\mu)$. Next,  by the above argument,  we have, for some subsequence
$(F_{n_k})_{k=1}^\infty $,
$(D_{xy}^{-+}g(F_{n_k}))(\gamma)\to (D_{xy} ^{-+}g(F))(\gamma)$ as $n\to\infty$ for
$\tilde\mu$-a.e.\ $(x,y,\gamma)$.

For any $x,y\in\R$, we evidently have
\begin{equation}\label{fuffu}
|g(x)-g(y)|\le|x-y|.\end{equation}
Therefore, the sequence $c(x,y,\gamma)^{1/2}(D^{-+}_{xy}
g(F_n))(\gamma)$, $n\in\N$, is $\tilde\mu$-uniformly
square-integrable, since so is the sequence
$c(x,y,\gamma)^{1/2}(D^{-+}_{xy}
F_n)(\gamma)$, $n\in\N$. Hence $$c(x,y,\gamma)^{1/2}(D_{xy}^{-+}g(F_{n_k}))(\gamma)\to c(x,y,\gamma)^{1/2}
(D_{xy}^{-+}g(F))(\gamma)\quad\text{as $k\to\infty$ in $L^2(\tilde\mu)$}.$$  By \eqref{iuedw} and  \eqref{guyiuhiu}, this yields: $g(F)\in D(\EK)$.

Finally, by  \eqref{iuedw}--\eqref{fuffu}, $\EK(g(F))\le\EK(F)$,
which means that $(\EK,D(\EK)))$ is a Dirichlet
form.\quad $\square$\vspace{2mm}

We shall now need  the bigger  space $\dd\Gamma$ consisting of all
$\Z_+\cup\{\infty\}$-valued Radon measures on $X$ (which is Polish, see e.g.\
\cite{Ka75}). Since $\Gamma\subset\dd\Gamma$ and ${\cal
B}(\dd\Gamma)\cap\Gamma={\cal B}(\Gamma)$, we can consider $\mu $
as a measure on $(\dd\Gamma,{\cal B}(\dd\Gamma))$ and
correspondingly $({\cal E}, D({\cal E}))$ as a Dirichlet form on
$L^2(\dd\Gamma,\mu)$.

For the notion of a quasi-regular Dirichlet form,  appearing in
the following lemma, we refer to \cite[Chap.~IV, Sect.~3]{MR}.

\begin{lem}\label{scfxgus}
 $(\EK,D(\EK))$ is a quasi-regular Dirichlet form on
$L^2(\dd\Gamma,\mu)$.
\end{lem}

\noindent {\it Proof}. Analogously to
\cite[Proposition~4.1]{MR98}, it suffices to show that there
exists a bounded, complete
 metric $\rho$ on $\dd\Gamma$ generating the vague
 topology such that, for all $\gamma_0\in\dd\Gamma$, $\rho(\cdot,\gamma_0)
\in D(\EK)$ and $$\int_X\gamma(dx)\int_X zm(dy)\, S(\rho(\cdot,\gamma_0))(x,y,\gamma)\le \eta(\gamma)\quad
\text{$\mu$-a.e.}$$  for some $\eta\in L^1(\dd\Gamma,\mu)$ (independent of
$\gamma_0$). Here, $S(F){:=}S(F,F)$. The proof below is a
modification of the proof of \cite[Proposition~4.8]{MR98} and the proof of \cite[Proposition~3.2]{KL}.

Fix any $x_0\in X$,   let $B(r)$ denote the open ball in $X$ of radius $r>0$ centered at $x_0$.
For each $k\in\N$, we define $$
g_k(x):=\,\frac23\left(\frac12-\operatorname{dist}(x,B(k))\wedge\frac12\right),\qquad
x\in X,$$ where $\operatorname{dist}(x,B(k))$ denotes the distance from the point
$x$ to the  ball $B(k)$.
Next, we set $$ \phi_k(x){:=} 3
g_k(x),\qquad x\in X,\ k\in\N.$$

Let $\zeta$ be a function in
$C_{\mathrm b}^1(\R)$ such that $0\le\zeta\le1$ on $[0,\infty)$,
$\zeta(t)=t$ on $[-1/2,1/2]$, $\zeta'\in[0,1]$ on $[0,\infty)$.
For any fixed $\gamma_0\in\dd\Gamma$ and for any $k,n\in\N$, (the
restriction to $\Gamma$ of) the function $$ \zeta\left(\sup_{j\le
n}|\la \phi_k g_j,\cdot\ra-\la\phi_kg_j,\gamma_0\ra|\right)$$
belongs to $\FC$ (note that $\la \phi_kg_j,\gamma_0\ra$ is a
constant).  Furthermore, taking into account that $\zeta'\in[0,1]$ on
$[0,\infty)$, we get from the mean value theorem, for each $
\gamma\in\Gamma$, $x\in\gamma$, and $x\in X\setminus \gamma$,
\begin{align}&S\left( \zeta\left(\sup_{j\le n}|\la \phi_k g_j,\cdot\ra-\la\phi_kg_j,\gamma_0\ra|\right)\right)(x,y,\gamma)\notag\\ &
\qquad \le c(x,y,\gamma)
\bigg( \sup_{j\le n}|\la \phi_k
g_j,\gamma\ra-\la\phi_kg_j,\gamma_0\ra -(\phi_kg_j)(x)+(\phi_kg_j)(y)|
\notag\\ &\qquad\quad-\sup_{j\le n}|\la \phi_k g_j,\gamma\ra-\la\phi_kg_j,\gamma_0\ra
|\bigg)^2 \notag \\ &\qquad \le c(x,y,\gamma) \sup_{j\le n}|-(\phi_kg_j)(x)+(\phi_kg_j)(y)|^2\notag\\
&\qquad\le 2c(x,y,\gamma)\bigg(\sup_{j\le n}(\phi_kg_j)(x)^2 +\sup_{j\le n}(\phi_kg_j)(y)^2 \bigg)\notag\\ &\qquad
\le 2 c(x,y,\gamma)(\pmb 1_{B(k+1/2)}(x)+\pmb 1_{B(k+1/2)}(y)).\label{joijhoi}
\end{align}

For each $k\in\N$, we define $$ F_k(\gamma,\gamma_0):=\zeta \left(\sup_{j\in\N}|\la \phi_kg_j,\gamma\ra-\la\phi_k g_j,\gamma_0\ra|\right),\qquad \gamma,\gamma_0\in\dd\Gamma.$$
Then, for a fixed $\gamma_0\in\dd\Gamma$, $$ \zeta \left(\sup_{j\le n}|\la \phi_kg_j,\gamma\ra-\la\phi_k g_j,\gamma_0\ra|\right)\to F_k(\gamma,\gamma_0) $$ as $n\to\infty$
for each $\gamma\in\dd\Gamma$ and in $L^2(\mu)$. Hence, by \eqref{joijhoi} and  the Banach--Alaoglu and  the Banach--Saks theorems (see e.g.\ \cite[Appendix~A.2]{MR}), $F_k(\cdot,\gamma_0)\in D(\EK)$ and $$ S(F_k(\cdot,\gamma_0))(x,y,\gamma)\le 2 c(x,y,\gamma) (\pmb 1_{B(k+1/2)}(x)+\pmb 1_{B(k+1/2)}(y))\quad\text{$\tilde
\mu$-a.e.}$$

Define $$ c_k:= \bigg(1+2\int c(x,y,\gamma)(\pmb 1_{B(k+1/2)}(x)+\pmb 1_{B(k+1/2)}(y))\, \tilde \mu(dx,dy,d\gamma)\bigg)^{-1/2}2^{-k/2},\qquad k\in\N,$$ which are finite positive numbers by \eqref{gyu}, and furthermore, $c_k\to0$ as $k\to\infty$.

 We define $$
\rho(\gamma_1,\gamma_2){:=}\sup_{k\in\N}
\big(c_kF_k(\gamma_1,\gamma_2)\big),\qquad
\gamma_1,\gamma_2\in\dd\Gamma.$$ By \cite[Theorem~3.6]{MR98},
$\rho$ is a bounded, complete metric on $\dd\Gamma$
generating the vague topology.

Analogously to the above, we now conclude that, for any fixed $\gamma_0\in\dd\Gamma$, $\rho(\cdot,\gamma_0)\in D(\EK)$
and $$ \int_X\gamma(dx)\int_X zm(dy) S(\rho(\cdot,\gamma_0))(x,y,\gamma)\le \eta(\gamma)\quad \text{$\mu$-a.e.,}$$
where $$ \eta(\gamma):= 2\sup_{k\in\N}\bigg(c_k^2\int_X\gamma(dx)\int_Xzm(dy)c(x,y,\gamma)(\pmb 1_{B(k+1/2)}(x)+\pmb 1_{B(k+1/2)}(y))\bigg).  $$
Finally, \begin{align*} \int_\Gamma \eta(\gamma)\,\mu(d\gamma)&\le 2\sum_{k=1}^\infty c_k^2\int c(x,y,\gamma)(\pmb 1_{B(k+1/2)}(x)+\pmb 1_{B(k+1/2)}(y))\,\tilde\mu(dx,dy,d\gamma)\\ &\le \sum_{k=1}^\infty 2^{-k}=1.\end{align*}
Thus, the lemma is proved.\quad $\square$\vspace{2mm}

For the notion of an exceptional set, appearing in
the next proposition, we refer  e.g.\  to \cite[Chap.~III,
Sect.~2]{MR}.

\begin{lem}\label{fdj}  The set $\dd\Gamma\setminus\Gamma$ is  $\EK$-exceptional\rom.
\end{lem}

\noindent {\it Proof}. We modify the proof of \cite[Proposition~1
and Corollary~1]{RS98} and the proof of \cite[Proposition~3.3]{KL} according to our situation.

It suffices to prove the result locally, i.e., to show that, for any fixed $a\in X$,
there exists a closed set $B_a$ that is the closure of an open neighborhood  of $a$ and
 such that the set $$
N_a:=\{\gamma\in\dd\Gamma: \sup_{x\in B_a}\gamma(\{x\})\ge2\}$$ is $\EK$-exceptional.
By
\cite[Lemma~1]{RS98}, we need to prove that there exists a
sequence $u_n\in D(\EK)$, $n\in\N$, such that each $u_n$ is a
continuous  function on $\dd\Gamma$, $u_n\to \pmb1 _{N_a}$ pointwise
as $n\to\infty$, and $\sup_{n\in\N}\EK (u_n)<\infty$.

So, we fix $a\in X$. There exists an open neighborhood $\tilde B_a$ of $a$ which is diffeomorphic
to the open cube $(-3,+3)^d$ in $\R^d$. We fix the corresponding coordinate system in $\tilde B_a$
and we set $B_a:=[-1,1]^d$.

Let $f\in C_0({\Bbb R})$ be such that $\pmb1_{[0,1]}\le f\le
\pmb1_{[-1/2,3/2)}$. For any $n\in\N$ and
$i=(i_1,\dots,i_d)\in{\cal A}_n:=\Z^d\cap[-n,n]^d$, we define a function $f_i^{(n)} \in C_0(X)$ by $$ f_i^{(n)}(x):=
\begin{cases}
\prod_{k=1}^d f(n x_k- i_k),&
  x\in \tilde B_a,\\ 0,&\text{otherwise}.\end{cases}$$ Let also $$I_i^{(n)}(x):=\begin{cases}\prod_{k=1}^d \pmb1
_{[-1/2,3/2)}(nx_k-i_k),& x\in\tilde B_a,\\0,&\text{otherwise},\end{cases}$$ and note that
$f_i^{(n)}\le I_i^{(n)}$.

Let $\psi\in C_{\mathrm b}^1 ({\Bbb R})$ be such that $\pmb1
_{[2,\infty)}\le\psi\le\pmb1_{[1,\infty)}$ and $0\le\psi'\le 2\,
\pmb1_{ (1,\infty)}$.
We define continuous functions $$\dd\Gamma\ni\gamma\mapsto
u_n(\gamma){:=}\psi \left( \sup_{i\in{\cal A}_n}\la
f_i^{(n)},\gamma\ra\right),\qquad n\in\N,$$ whose restriction to
$\Gamma$ belongs to $\FC$. Evidently, $u_n\to \pmb1_{N_a}$ pointwise as
$n\to\infty$.

By the mean value theorem, we have, for each $\gamma\in\Gamma$, $x\in\gamma$, $y\in X\setminus
\gamma$, and for some point $T_n(x,y,\gamma)$ between $\sup_{i\in{\cal A}_n}\la f_i^{(n)},\gamma\setminus x\cup y\ra$ and $\sup_{i\in{\cal A}_n}\la f_i^{(n)},\gamma\ra$:
\begin{gather}
S(u_n)(x,y,\gamma)=c(x,y,\gamma)\psi'(T_n(x,y,\gamma))^2\bigg(\sup_{i\in{\cal A}_n}\la f_i^{(n)},\gamma\setminus x\cup y\ra- \sup_{i\in{\cal A}_n}\la f_i^{(n)},\gamma\ra\bigg)^2\notag\\
\le c(x,y,\gamma)\psi'(T_n(x,y,
\gamma))^2 \sup_{i\in{\cal A}_n}|\la f_i^{(n)},\gamma\setminus x\cup y\ra- \la f_i^{(n)},\gamma\ra|^2\notag\\
\le 2c(x,y,\gamma)\psi'(T_n(x,y,\gamma))^2\bigg(\sup_{i\in{\cal A}_n}f_i^{(n)}(x)^2+\sup_{i\in{\cal A}_n}f_i^{(n)}(y)^2\bigg)\notag\\ \le 8c(x,y,\gamma)  (\pmb1_{\tilde B_a}(x)+\pmb1_{\tilde B_a}(y)). \label{yuftyu}
\end{gather}
By \eqref{gyu}  and \eqref{yuftyu}, we  conclude that $$ \sup_{n\in\N}\EK(u_n)< \infty,$$ which implies the lemma.\quad $\square$

We now have the main result of this paper.

 \begin{theorem}\label{8435476}
 There
exists a conservative Hunt process 
$${\bf M}=({\pmb{ \Omega}},{\bf F},({\bf F}_t)_{t\ge0},({\pmb
\Theta}_t)_{t\ge0}, ({\bf X}(t))_{t\ge 0},({\bf P
}_\gamma)_{\gamma\in\Gamma})$$ on $\Gamma$ \rom(see e\rom.g\rom.\
\rom{\cite[p.~92]{MR})} which is properly associated with $({\cal
E},D({\cal E}))$\rom, i\rom.e\rom{.,} for all \rom($\mu$-versions
of\/\rom) $F\in L^2(\Gamma,\mu)$ and all $t>0$ the function
\begin{equation}\label{zrd9665} \Gamma\ni\gamma\mapsto
p_tF(\gamma){:=}\int_{\pmb\Omega} F({\bf X}(t))\, d{\bf
P}_\gamma\end{equation} is an ${\cal E}$-quasi-continuous version
of $\exp(-t{H})F$\rom, where $(H,D(H))$ is the generator of $({\cal
E},D({\cal E}))$\rom.  $\bf M$ is up to $\mu$-equivalence unique
\rom(cf\rom.\ \rom{\cite[Chap.~IV, Sect.~6]{MR}).} In
particular\rom, ${\bf M}$ is $\mu$-symmetric \rom(i\rom.e\rom{.,}
$\int G\, p_tF\, d\mu=\int F \, p_t G\, d\mu$ for all
$F,G:\Gamma\to{\Bbb R}_+$\rom, ${\cal B}(\Gamma)$-measurable\rom)\rom, so
 has $\mu$ as an invariant measure\rom.

\rom{2)} ${\bf M}$ from \rom{1)} is  up to $\mu$-equivalence
\rom(cf\rom.\ \rom{\cite[Definition~6.3]{MR}}\rom) unique between
all Hunt processes ${\bf M}'=({\pmb{ \Omega}}',{\bf F}',({\bf
F}'_t)_{t\ge0},({\pmb \Theta}'_t)_{t\ge0}, ({\bf X}'(t))_{t\ge
0},({\bf P }'_\gamma)_{\gamma\in\Gamma})$ on $\Gamma$ having $\mu$
as  invariant measure and solving the martingale problem for
$(-H, D(H))$\rom, i\rom.e\rom.\rom, for all $G\in D(H)$
$$\widetilde G({\bf X}'(t))-\widetilde G({\bf X}'(0))+\int_0^t (H
G)({\bf X}'(s))\,ds,\qquad t\ge0,$$ is an $({\bf
F}_t')$-martingale under ${\bf P}_\gamma'$ for ${\cal
E}$-q\rom.e\rom.\ $\gamma\in\Gamma$\rom. \rom(Here\rom,
$\widetilde G$ denotes an ${\cal E}$-quasi-continuous version of $G$\rom,
cf\rom. \rom{\cite[Ch.~IV, Proposition~3.3]{MR}.)}
\end{theorem}

\begin{rem}\rom{
 In Theorem~\ref{8435476}, ${\bf M}$ can be taken canonical, i.e., $\pmb\Omega$ is
the set of all {\it cadlag} functions $\omega:[0,\infty)\to
\Gamma$ (i.e., $\omega$ is right continuous on $[0,\infty)$ and
has left limits on $(0,\infty)$), ${\bf
X}(t)(\omega){:=}\omega(t)$, $t\ge 0$, $\omega\in\pmb\Omega$,
$({\bf F}_t)_{t\ge 0}$ together with ${\bf F}$ is the corresponding
minimum completed admissible family (cf.\
\cite[Section~4.1]{Fu80}) and ${\pmb \Theta}_t$, $t\ge0$, are the
corresponding natural time shifts.
 }\end{rem}

\noindent {\it Proof of Theorem\/}~\ref{8435476}. The first part
of the theorem follows from
 Lemmas~\ref{guzazagus}--\ref{fdj}, the fact that $1\in D({\cal E})$,  ${\cal E}(1,1)=0$,  and
\cite[Chap.~IV, Theorem~3.5 and Chap.~V, Proposition~2.15]{MR}.
The second part follows directly from (the proof of)
\cite[Theorem~3.5]{AR}.\quad $\square$
\vspace{2mm}

Let us now derive an explicit formula for the generator of  $\EK$. However, this can only be done under stronger conditions on the coefficient $c(x,y,\gamma)$.

 Using \eqref{fdrtsdrt} and \eqref{jkhgiu}, we have, for $F\in\FC$,
\begin{gather}\EK(F) =\int_\Gamma \mu(d\gamma)\int_X zm(dy)\exp[-E(y,\gamma)+E(y,\gamma)]
 \int_X\gamma(dx)\,
c(x,y,\gamma) (D_{xy}^{-+}F)(\gamma)^2\notag\\ =\int_\Gamma \mu(d\gamma)\int_X\gamma(dy)
\exp[E(y,\gamma\setminus y)]
\int_X (\gamma\setminus y)(dx)c(x,y,\gamma\setminus y)(F(\gamma\setminus x)-F(\gamma\setminus y))^2\notag\\
=\int_\Gamma\mu(d\gamma)\int_X\gamma(dx)\int_X (\gamma\setminus x)(dy)\exp[E(y,\gamma\setminus y)]
 c(x,y,\gamma\setminus y)(F(\gamma\setminus x)-F(\gamma\setminus y))^2\notag\\ =\int_\Gamma
\mu(d\gamma)\int_X zm(dx)\exp[-E(x,\gamma)]\int_X\gamma(dy)\exp[E(y,\gamma\setminus y\cup x)]\notag\\ \times
c(x,y,\gamma\setminus y\cup x)(F(\gamma)-F(\gamma\setminus y\cup x))^2
\notag\\ =\int_\Gamma
\mu(d\gamma)\int_X\gamma(dx)\int_X zm(dy)\,c(y,x,\gamma\setminus x\cup y)\notag \\ \times \exp[-E(y,\gamma)+ E(x,\gamma\setminus x\cup y)] (D_{xy}^{-+}F)(\gamma)^2.\label{yfu}
\end{gather}
By \eqref{2} and \eqref{yfu}, we have, for any $F,G\in\FC$,
$$\EK(F,G)= \int_\Gamma
\mu(d\gamma)\int_X\gamma(dx)\int_X zm(dy) \, \tilde c(x,y,\gamma) (D_{xy}^{-+}F)(\gamma)(D_{xy}^{-+}G)(\gamma),$$
where \begin{equation}\label{kgyuytty}\tilde c(x,y,\gamma)=\frac12\big(
c(x,y,\gamma)+c(y,x,\gamma\setminus x\cup y)\exp[-E(y,\gamma)+
 E(x,\gamma\setminus x\cup y)]\big). \end{equation}
As easily seen, $\tilde c$
again satisfies the condition \eqref{gyu}. Furthermore, $\tilde c$ evidently satisfies the following identity:
$$ \tilde c(x,y,\gamma)=\tilde c(y,x,\gamma\setminus x\cup y)\exp[-E(y,\gamma)+
 E(x,\gamma\setminus x\cup y)],$$
so that $\tilde{\tilde c}=\tilde c$.

\begin{theorem} \label{jihgy}
 Assume that\rom, for each compact $\Lambda\subset X$\rom, \begin{equation}\label{jkg}
\int_X\gamma(dx)\int _X zm(dy)\, \tilde c(x,y,\gamma)(\pmb1_\Lambda(x)+\pmb1_\Lambda(y))\in L^2(\Gamma,\mu).
\end{equation}
 Then
\begin{equation}\label{ghjfytr} \EK(F,K)=\int_\Gamma
(HF)(\gamma)G(\gamma)\, \mu(d\gamma),\qquad F,G\in\FC,\end{equation}
where
\begin{equation}\label{eghdn}
(HF)(\gamma)=-2\int_{X}\gamma(dx)\int _X zm(dy) \tilde c(x,y,\gamma) (D^{-+}_{xy}F)(\gamma)\qquad \text{\rom{$\mu$-a.e.}}
\end{equation}
and $HF\in L^2(\Gamma,\mu)$\rom.  The Friedrichs' extension of the operator $(H,\FC)$
in $L^2(\Gamma,\mu)$ is  $(H,D(H))$\rom.

\end{theorem}

\noindent {\it Proof.} Formulas  \eqref{ghjfytr} and \eqref{eghdn} follow from
\eqref{fdrtsdrt} and  \eqref{kgyuytty}, analogously to \eqref{yfu}.
The fact that  $HF\in L^2(\Gamma,\mu)$ trivially follows from  \eqref{jkg}. \quad $\square$

\section{Examples}\label{lkhuhy}

Throughout this section, we shall always assume that a pair potential $\phi$, an activity $z$, and a corresponding Gibbs measure $\mu\in {\cal G}(z,E)$ are  either as in Theorem~\ref{gdtrsresa}
 or as in Theorem~\ref{yte6te645}. Furthermore, in the case $X=\R^d$, we shall also suppose that the condition of Lemma~\ref{lohkjj} is satisfied. Thus, in any case we have that  $\phi$ is bounded from below,  satisfies (I), and $$\sum_{y\in\gamma}|\phi(x,y)|<\infty\quad\text{for $m\otimes \mu$-a.e.\ $(x,\gamma)\in X\times\Gamma$}.$$
 By \eqref{fdrtsdrt}, the latter easily implies that, for $\mu$-a.e.\ $\gamma\in\Gamma$ and for each $x\in\gamma$,
 $$\sum_{y\in\gamma\setminus x}|\phi(x,y)|<\infty.$$

We shall now consider some examples of the coefficient $c(x,y,\gamma)$ for which the above assumptions are satisfied.

Let $a:X^2\to[0,\infty]$ be a  symmetric measurable function such that \begin{equation}\label{diguy}
\sup_{x\in\Lambda}\int_X a(x,y)\,m(dy)<\infty, \quad \sup_{y\in\Lambda}\int_X a(x,y)\,m(dx)<\infty\end{equation}
for any compact $\Lambda\subset X$.

\begin{rem}
\rom{ In the case $X=\R^d$, it is natural to suppose that the function $a$ is translation invariant, i.e., $a(x,y)=\tilde a(x-y)$
for some $\tilde a:X\to[0,\infty]$, $\tilde a(-x)=\tilde (x)$, $x\in\R^d$, in which case \eqref{diguy} is equivalent to the integrability of $\tilde a$.
}\end{rem}

For $s\in[0,1]$
we define \begin{equation}\label{ufytre} c(x,y,\gamma)=c_s(x,y,\gamma):=a(x,y)\exp[sE(x,\gamma\setminus x)-(1-s)E(y,\gamma\setminus x)]
.\end{equation}
We evidently have  $\tilde c_s (x,y,\gamma) =c_s(x,y,\gamma)$.

\begin{prop} \label{hgiuyg} \rom{1)} For each $s\in [0,1]$\rom, the coefficient $c_s$ satisfies \eqref{gyu}\rom.

\rom{2)} Assume  that the function $a$ is bounded\rom. Then\rom, for each $s\in[0,1/2]$\rom,
the coefficient $c_s$ satisfies \eqref{jkg}\rom.
Furthermore\rom, for each $s\in (1/2,1]$\rom,  \eqref{jkg} is satisfied if additionally \begin{equation}\label{hgft}
\sup_{x\in X}\int_X|\exp[(2s-1)\phi(x,y)]-1|\,m(dy)<\infty.
\end{equation}

\end{prop}

\noindent {\it Proof.}
1) We have \begin{align}
&\int_\Gamma \mu(d\gamma)\int_X\gamma(dx)\int_Xm(dy) a(x,y)\notag\\
&\qquad\times \exp[sE(x,\gamma\setminus x)-(1-s)E(y,\gamma\setminus x)](\pmb1_\Lambda(x)+\pmb1_\Lambda(y))\notag\\
&=\int_\Gamma \mu(d\gamma)\int_X zm(dx)\int_Xm(dy) a(x,y)
\notag\\ &\quad\times\exp[(s-1)E(x,\gamma)+(s-1)E(y,\gamma)] )](\pmb1_\Lambda(x)+\pmb1_\Lambda(y))\notag\\
&=\int_X zm(dx)\int_Xm(dy) a(x,y)(\pmb1_\Lambda(x)+\pmb1_\Lambda(y))\notag\\
&\quad\times \int_\Gamma \mu(d\gamma)\prod_{u\in\gamma}
(1+(\exp[(s-1)\phi(x,u)+(s-1)\phi(y,u)]-1))\notag\\
&=\int_X zm(dx)\int_Xm(dy) a(x,y)(\pmb1_\Lambda(x)+\pmb1_\Lambda(y))\notag\\
&\quad\times \int_\Gamma \mu(d\gamma)\bigg( 1+\sum_{n=1}^\infty
\sum_{\{u_1,\dots,u_n\}\subset\gamma}\prod_{i=1}^n (\exp[(s-1)\phi(x,u_i)+(s-1)\phi(y,u_i)]-1)\bigg)\notag\\
&=\int_X zm(dx)\int_Xm(dy) a(x,y)(\pmb1_\Lambda(x)+\pmb1_\Lambda(y))\notag\\
&\quad\times \int_\Gamma \mu(d\gamma)\bigg( 1+\sum_{n=1}^\infty
\frac1{n!}\int_{X^n}\prod_{i=1}^n (\exp[(s-1)\phi(x,u_i)+(s-1)\phi(y,u_i)]-1) \notag\\
&\quad\times k_\mu^{(n)}(u_1,\dots,u_n)\,m(du_1)\dotsm m(du_n)\bigg).\label{1234}
\end{align}
Using the Ruelle bound, we get, for any $x,y\in X$, 
\begin{align*}
&\int_{X^n}\prod_{i=1}^n |\exp[(s-1)\phi(x,u_i)+(s-1)\phi(y,u_i)]-1|\\
&\qquad \times k_\mu^{(n)}(u_1,\dots,u_n)\,m(du_1)\dotsm m(du_n)\\
&\quad\le \bigg(\xi \int_{X} |\exp[(s-1)\phi(x,u)+(s-1)\phi(y,u)]-1|\,m(du)\bigg)^n \\
&\quad\le \bigg(\xi \int_{X} |\exp[-\phi(x,u)-\phi(y,u)]-1|\,m(du)\bigg)^n \\
&\quad \le\bigg( \xi \int_X\big( |\exp[-\phi(x,u)]-1|+\exp[-\phi(x,u)]\,
|\exp[-\phi(y,u)]-1|\big)m(du)\bigg)^n.
\end{align*} 
From here, (I) and \eqref{1234}, the statement follows.

2) Analogously,   we have:
\begin{gather}
\int_\Gamma \mu(d\gamma)\bigg(\int_X\gamma(dx)\int_X zm(dy)c_s(x,y,\gamma)(\pmb1_\Lambda(x)+\pmb1_\Lambda(y))\bigg)^2\notag\\
=\int_\Gamma\mu(d\gamma)\int_X \gamma(dx)\int_X zm(dy)\int_X zm(dy')c_s(x,y,\gamma)c_s(x,y',\gamma)\notag\\ \times
(\pmb1_\Lambda(x)+\pmb1_\Lambda(y))(\pmb1_\Lambda(x)+\pmb1_\Lambda(y'))\notag\\
\text{}+\int_\Gamma \mu(d\gamma)\int_X\gamma(dx)\int_X (\gamma\setminus x)(dx')\int_X zm(dy)
\int_X zm(dy')\notag\\ \times
c_s(x,y,\gamma)c_s(x',y',\gamma)(\pmb1_\Lambda(x)+\pmb1_\Lambda(y))
(\pmb1_\Lambda(x')+\pmb1_\Lambda(y'))\notag\\
=\int_\Gamma \mu(d\gamma)\int_X zm(dx)\int_X zm(dy)\int_X zm(dy')\exp[-E(x,\gamma)]\notag\\
\times c_s(x,y,\gamma\cup x)c_s(x,y',\gamma\cup x)
(\pmb1_\Lambda(x)+\pmb1_\Lambda(y))(\pmb1_\Lambda(x)+\pmb1_\Lambda(y'))\notag\\
\text{}+ \int_\Gamma\mu(d\gamma)\int_X zm(dx) \int_X zm(dx') \int_X zm(dy) \int_X zm(dy')\notag\\ \times
\exp[-E(x,\gamma)-E(x',\gamma)-\phi(x,x')] \notag\\ \times c_s(x,y,\gamma\cup x\cup x')  c_s
(x',y',\gamma\cup x\cup x')(\pmb1_\Lambda(x)+\pmb1_\Lambda(y))(\pmb1_\Lambda(x')+\pmb1_\Lambda(y'))\notag\\
=\int_X zm(dx)\int_X zm(dy)\int_X zm(dy') a(x,y) a(x,y') (\pmb1_\Lambda(x)+\pmb1_\Lambda(y))
(\pmb1_\Lambda(x)+\pmb1_\Lambda(y')) \notag\\ 
\times \int_\Gamma\mu(d\gamma)\exp\bigg[\sum_{u\in\gamma}\big(
(2s-1)\phi(x,u)-(1-s)\phi(y,u)-(1-s)\phi(y',u)
\big)\bigg]\notag\\ \text{}+\int_X zm(dx) \int_X zm(dx') \int_X zm(dy) \int_X zm(dy') a(x,y)a(x',y')\notag\\ \times
(\pmb1_\Lambda(x)+\pmb1_\Lambda(y))(\pmb1_\Lambda(x')+\pmb1_\Lambda(y'))\notag\\ \times
\exp[(2s-1)\phi(x,x')
-(1-s)\phi(x',y) -(1-s)\phi(x,y')
]
\notag\\ \times \int_\Gamma\mu(d\gamma)\exp\bigg[\sum_{u\in\gamma}-(1-s)(
\phi(x,u)+\phi(x',u)+\phi(y,u)+\phi(y',u))\bigg]\notag\\ \le
C_2\bigg(\int_X zm(dx)\int_X zm(dy)\int_X zm(dy') a(x,y) a(x,y') (\pmb1_\Lambda(x)+\pmb1_\Lambda(y))
(\pmb1_\Lambda(x)+\pmb1_\Lambda(y')) \notag\\
\times \exp\bigg[\xi\sup_{x\in X}\sup_{y\in X}\sup_{y'\in X}  \int_X |\exp[(2s-1)\phi(x,u)\notag\\ \text{}-(1-s)\phi(y,u)-(1-s)\phi(y',u)]-1|\,m(du)\bigg]\notag\\
\text{}+\int_X zm(dx) \int_X zm(dx') \int_X zm(dy) \int_X zm(dy') a(x,y)a(x',y')\notag\\ \times
(\pmb1_\Lambda(x)+\pmb1_\Lambda(y))(\pmb1_\Lambda(x')+\pmb1_\Lambda(y'))\exp[(2s-1)\phi(x,x')]
\notag\\ \times \exp\bigg[\xi\sup_{x\in X}\sup_{x'
\in X}\sup_{y\in X}\sup_{y'\in X}\int_X|\exp[-(1-s)(
\phi(x,u)+\phi(x',u)\notag\\ \text{}+\phi(y,u)+\phi(y',u))]-1|\,m(du)\bigg]\bigg)\notag\\ \le C_3\bigg(
\int_X zm(dx) \int_X zm(dy)\int_X zm(dy') a(x,y) a(x,y') (\pmb1_\Lambda(x)+\pmb1_\Lambda(y))
(\pmb1_\Lambda(x)+\pmb1_\Lambda(y')) \notag\\ \text{}+\int_X zm(dx) \int_X zm(dx') \int_X zm(dy) \int_X zm(dy') a(x,y)a(x',y')\notag\\ \times
(\pmb1_\Lambda(x)+\pmb1_\Lambda(y))(\pmb1_\Lambda(x')+\pmb1_\Lambda(y'))|\exp[(2s-1)\phi(x,x')]-1|\notag\\
\text{}+ \bigg(\int_Xzm(dx)\int_Xzm(dy)a(x,y)(\pmb1_\Lambda(x)+\pmb1_\Lambda(y))\bigg)^2
\bigg),\label{iugtut} \end{gather} where $C_2,C_3>0$.
Using (I), \eqref{diguy}, \eqref{hgft}, and the boundedness of $a$, we easily conclude that the expression in \eqref{iugtut}
is finite. Indeed, for example, we have:
\begin{align*}&\int_X zm(dx)\int_X zm(dx')\int_\Lambda zm(dy)\int_\Lambda zm(dy')a(x,y)a(x',y')|\exp[(2s-1)\phi(x,x')]-1|\\
&\qquad\le \bigg(\sup_{(u,v)\in X^2}a(u,v)\bigg) \int_\Lambda zm(dy)\int_\Lambda zm(dy')\\ &\qquad\qquad \times \int_Xzm(dx)a(x,y)\int_X zm(dx')|\exp[(2s-1)\phi(x,x')]-1|<\infty.
\end{align*}
Thus, the proposition is proved.\quad $\square$\vspace{2mm}

Let us now present a straightforward generalization of the above result. Let now $a:X^2\to\R$ be a measurable function  which satisfies \eqref{diguy} (and which is not necessarily symmetric). For $u,v\in[0,1]$, we define
$$ \varkappa(x,y,\gamma)= \varkappa_{u,v}(x,y,\gamma):=
\exp[uE(x,\gamma\setminus x)-(1-v)E(y,\gamma)],$$ and $$ c(x,y,\gamma)=c_{u,v}(x,y,\gamma):=a(x,y)\varkappa_{u,v}(x,y,\gamma)
.$$ In particular, for $u=v$, we get the previous example of a Kawasaki dynamics. Note also that, for $u=0$ and $v=1$, we get $$c_{0,1}(x,y,\gamma)=a(x,y). $$

By \eqref{kgyuytty}, we have \begin{align*}
\tilde c_{u,v}(x,y,\gamma)&=\frac12\big(a(x,y )\exp[uE(x,\gamma\setminus x)-(1-v)E(y,\gamma)]\\ &\qquad+a(y,x)
\exp[vE(x,\gamma\setminus x)-(1-u)E(y,\gamma)]\big).
\end{align*}

Absolutely analogously to Proposition~\ref{hgiuyg}, one can prove its following generalization.

\begin{prop}\rom{1)} For each $u,v\in [0,1]$\rom, the coefficient $c_{u,v}$ satisfies \eqref{gyu}\rom. 

\rom{2)} Assume  that the function $a$ is bounded\rom. Then\rom,
\eqref{jkg} is satisfied if $$
\sup_{x\in X}\int_X|\exp[(2(u\vee v)-1)\phi(x,y)]-1|\,m(dy)<\infty.
$$
\end{prop}

\section{Scaling limits of Kawasaki dynamics}\label{newsection}

We start this section with a brief discussion of Glauber dynamics of continuous particle systems. 

\subsection{Glauber (birth-and-death) dynamics}

In the classical  Ising model, the Glauber dynamics means that particles randomly change their spin value, which is called a spin-flip. The generator of this dynamics is given by
$$ (H_{\mathrm G}f)(\sigma)=\sum_{x\in\Z^d}a(x,\sigma)(\nabla_xf)(\sigma),$$
where $$ (\nabla_xf)(\sigma)=f(\sigma^x)-f(\sigma),$$
$\sigma^x$ denoting the configuration $\sigma$ in which the particle at site $x$ has  changed its spin value. In the interpretation of a lattice system with spin space $S=\{-1,1\}$ as a model of a lattice gas, 
the Glauber dynamics means that, at each site $x$, a particle randomly appears and disappears. Hence, this dynamics may be interpreted as a birth-and-death process on $\Z^d$. Therefore, in the continuous case, an analog of  the Glauber dynamics should be a process in which particles randomly appear and disappear in the space, i.e., 
a spatial birth-and-death process. The generator of such a process is informally given by the formula $$
(H_{\mathrm G}F)(\gamma)=\sum_{x\in\gamma}d(x,\gamma)(D^-_xF)(\gamma)+\int _{\R^d}b(x,\gamma) (D^+_xF)(\gamma)
\,dx,$$ 
 where $$ (D_x^-F)(\gamma)=F(\gamma\setminus x)-F(\gamma),\quad 
 (D_x^+F)(\gamma)=F(\gamma\cup x)-F(\gamma).$$

Spatial birth-and-death processes were first discussed by Preston in \cite{P}. Under some conditions on the birth and death rates, Preston proved the existence  of such  processes in a bounded domain in $\R^d$. Though the number of particles can be arbitrarily large in this case, the total number of particles remains finite at any moment of time.  

The
problem of construction of a spatial birth-and-death process in
the infinite volume was initiated by  Holley and  Stroock in  \cite{HS}.
In fact, in that paper, birth-and-death processes in bounded domains were analyzed in detail. Only in a very special case of nearest neighbor birth-and-death processes on the real line, the existence of a corresponding process on the whole space was proved and its properties were studied. 
In \cite{G1}, Gl\"otzl  derived conditions on the coefficients $d(x,\gamma)$, $b(x,\gamma)$ under which the Glauber  generator becomes a symmetric operator in the space $L^2(\mu)$, where $\mu$ is a given Gibbs measure. Let us also mention the papers \cite{BCC,Wu} devoted to the study of the spectral gap of the Glauber dynamics in the finite volume, for which the death coefficient is equal to 1. An analog of such a  dynamics, but on the whole space (thus, involving infinite configurations), was  constructed in \cite{KL}.  The coefficients of the generator of this dynamics are given by $$ d(x,\gamma)=1,\quad b(x,\gamma)=\exp[-E(y,\gamma)] ,$$
and this dynamics has a Gibbs measure corresponding to the pair potential $\phi$ as symmetrizing measure. 
The result
about the spectral gap for a positive $\phi$  has also been extended in \cite{KL} to the infinite volume. 
We also refer to \cite{Martin} for a discussion of a scaling limit of equilibrium fluctuations of this  dynamics.

By analogy with the Kawasaki dynamics, we are now able to construct 
an equilibrium Kawasaki dynamics in the general case. So, we consider
a measurable mapping $$ X\times\Gamma\ni(x,\gamma)\mapsto d(x,\gamma)\in[0,\infty)$$ and assume that, for each compact $\Lambda\subset X$,
\begin{equation}\label{uftfr} \int_\Gamma \mu(d\gamma)\int_\Lambda \gamma(dx)d(x,\gamma)<\infty.\end{equation}
We define a bilinear form
$$ \EG(F,G):= \int_\Gamma\mu(d\gamma)\int_X \gamma(dx) d(x,\gamma) (D_x^-F)(\gamma) (D_x^-G)(\gamma),$$
where $F,G\in\FC$. This bilinear form is closable on $L^2(\Gamma,\mu)$, and its closure will be denoted by $(\EG,D(\EG))$. The latter is a quasi-regular Dirichlet form on $L^2(\dd\Gamma,\mu)$. Furthermore, the set $\dd\Gamma\setminus\Gamma$ is $\EG$-exceptional.
Therefore,  there
exists a conservative Hunt process 
which is properly associated with $({\cal
E}_{\mathrm G},D({\cal E}_{\mathrm G}))$.  

Next, assume that, for each compact $\Lambda\subset X$\rom, \begin{align}
\int_\Lambda \gamma(dx) d(x,\gamma)&\in L^2(\Gamma,\mu),\notag\\ 
\int_\Lambda zm(dx)b(x,\gamma)&\in L^2(\Gamma,\mu),\label{asd}\end{align}
where $$ b(x,\gamma):= \exp[-E(x,\gamma)]d(x,\gamma\cup x),\qquad x\in X, \gamma\in\Gamma.$$ Then, $\FC$ is a subset of the domain of the generator $H_{\mathrm G}$ of the Dirichlet form $(\EG,D(\EG))$, and for each $F\in\FC$, 
$$ (H_{\mathrm G}F)(\gamma)=-\int_{X}z
m(dx)\, b(x,\gamma) (D^+_xF)(\gamma)-\int_{X}\gamma(dx)\,d(x,\gamma)(D^-_xF)
(\gamma)\qquad \text{\rom{$\mu$-a.e.}}$$

Also by analogy with the Kawasaki dynamics, one  can construct the following examples of the Glauber dynamics (the Gibbs measure $\mu$ being the same as in Section~\ref{lkhuhy}). For each $s\in [0,1]$, we define $$ d(x,\gamma)=d_s(x,\gamma):= \exp[s E(x, \gamma\setminus x)],$$
so that $$ b(x,\gamma)=b_s(x,\gamma)=\exp[(s-1)E(x, \gamma)].$$
 Then, for each $s\in[0,1]$, \eqref{uftfr} holds, and therefore  the corresponding Glauber dynamics exists. Furthermore, for each $s\in[0,1/2]$, the coefficients $d_s$, $b_s$ satisfy \eqref{asd}, while for $s\in(1/2,1]$, \eqref{asd} is satisfied if 
\eqref{hgft} holds.

Note that, though the construction of the Glauber dynamics and that of the Kawasaki dynamics look quite similar, there is a drastic difference between them in that (at least heuristically) the law of conservation of the number of particles holds for the Kawasaki dynamics, and does not for the Glauber dynamics. We, therefore, cannot expect a spectral gap for the generator of the  Kawasaki dynamics in the infinite volume. 

Furthermore, the Glauber and Kawasaki dynamics have different sets 
 of symmetrizing measures. Indeed, the set of symmetrizing measures of a given Glauber dynamics consists of all grand-canonical Gibbs measures corresponding to a given relative energy of interaction and a {\it fixed} activity parameter $z>0$,
 while for a given Kawasaki dynamics activity parameter $z>0$ may be {\it arbitrary}. 
 
\subsection{Glauber dynamics as a limiting Kawasaki dynamics }
\label{yugfu}

Let $X=\R^d$ and let $\mu$ be a Gibbs measure as in Theorem~\ref{yte6te645}, 1) (low activity-high temperature regime). We fix a function $\tilde a:\R^d\to[0,\infty)$ such that $\tilde a(-x)=\tilde a(x)$, $x\in\R^d$, and $\tilde a\in L^1(\R^d,dx)$.  For each $s\in[0,1]$, consider the Kawasaki dynamics corresponding to the coefficient $c=c_s$
 given by \eqref{ufytre} with $a(x,y):=\tilde a(x-y)$. 
 
 Let us now consider the following scaling of this dynamics. For each $\delta>0$, define
 $$ \tilde a_\delta(x)=\delta^d\tilde a(\delta\cdot), \quad x\in\R^d,$$
 and consider the $(\delta,s)$-Kawasaki dynamics which is defined just as the dynamics above, but by using the function $\tilde a_\delta$ instead of $\tilde a$. Denote by $(H_{\delta,s},D(H_{\delta,s}))$ the generator of this dynamics.  
 
  Let us also fix the $s$-Glauber dynamics  
 corresponding to 
 \begin{align*}
d(x,\gamma)&=d_s(x,\gamma):=\alpha \exp[sE(x,\gamma\setminus x)] ,\\
b(x,\gamma)&=b_s(x,\gamma)=\alpha\exp[(s-1)E(x,\gamma)] ,\end{align*}
where $s\in[0,1]$ and 
$$ \alpha:= 2k_\mu^{(1)}\int_{\R^d}a(x)\,dx$$
(note that the first correlation function $k_\mu^{(1)}$ is a constant).  
Denote by $(H_{0,s},D(H_{0,s}))$ the generator of this dynamics.

 We expect that the $s$-Glauber dynamics is the limit of the $(\delta,s)$-Kawasaki dynamics as $\delta\to0$. 
  In particular, in the case $s=0$, it is shown in \cite{FKL} that, for each $\varphi\in C_0(\R^d)$,
   $$H_{\delta,s}e^{\langle\varphi,\cdot\rangle}\to H_{0,s}e^{\langle\varphi,\cdot\rangle}\quad \text{in $L^2(\Gamma,\mu)$} $$ as $\delta\to0$.  
  In the case where the potential  $\phi$ is non-negative, 
  one can conclude from \cite{KL} that the set of finite linear combinations of the exponential functions is a core for the Glauber generator $(H_{0,s},D(H_{0,s}))$. From here, using a classical result from the theory of semigroups \cite{Dav}, one derives the weak convergence of finite-dimensional distributions of the corresponding 
equilibrium dynamics, starting with their equilibrium  distribution  $\mu$.   
   
   \subsection{Diffusion approximation for  the Kawasaki dynamics}
 Now, let $X=\R^d$ and let $\mu$ be a Gibbs measure as in Theorem~\ref{yte6te645}. We will consider a Kawasaki dynamics as in subsec.~\ref{yugfu}, but this time we will additionally assume that 
 $\tilde a(x)$ only depends on $|x|$, and has compact support. We again consider the corresponding $(\delta,s)$-Kawasaki dynamics, but this time we are interested in its limiting behavior as $\delta\to\infty$. It appears that, we additionally have to re-scale time by multiplying it  by $\delta^2$. Thus, the generator of this dynamics is given by
 $ \widetilde H_{\delta,s}=\delta^2 H_{\delta,s}$. Under quite weak assumptions on the potential $\phi$, it is shown in \cite{KKL} that, for each function $F$ from some set of smooth local functions on the configuration space, 
 $$ \widetilde H_{\delta,s}F\to \widetilde H_sF\quad \text{in } L^2(\Gamma,\mu)$$ as $\delta\to\infty$. Here $\widetilde H_s$
 is the generator of a diffusion dynamics given by 
 \begin{align*}
 \widetilde H_{s} F(\gamma)&=c\int_{\R^d}\gamma(dx)\bigg(-\Delta_xF(\gamma)+2\sum_{u\in\gamma\setminus x}\langle
 \nabla_x F(\gamma),s\nabla\phi(x-u)\rangle\bigg)
 \\&\quad\times 
 \exp[(-2s+1)E(x,\gamma\setminus x)],
 \end{align*}  
 where
 $$ c:=z\int_{\R^d} a(x)(x^1)^2\,dx$$
 ($x^1$ denoting the first coordinate of $x\in\R^d$), $\Delta_x F(\gamma):=\Delta_yF(\gamma\setminus x\cup y)\big|_{y=x}$, and 
 $\nabla_xF(\gamma)$ is defined by analogy.  In particular, for $s=1/2$,
 $$ \widetilde H_{s} F(\gamma)=c\int_{\R^d}\gamma(dx)\bigg(-\Delta_xF(\gamma)+\sum_{u\in\gamma\setminus x}\langle
 \nabla_x F(\gamma),\nabla\phi(x-u)\rangle\bigg),$$
 which is the generator of the gradient stochastic dynamics, e.g.\ \cite{111,Fritz}.  
 
In the case where $\mu$ is a Gibbs measure as in Theorem~\ref{yte6te645}, 3), 
$\tilde\phi\in C_{\mathrm b}^3(\R^d)$, and $\tilde \phi$ sufficiently quickly converges to zero at infinity, Choi,  Park, and  Yoo \cite{CPY}
found a core for the generator of the gradient stochastic dynamics.  
Using this result, we  derive the weak convergence of finite-dimensional distributions of the corresponding 
equilibrium dynamics, starting with their equilibrium  distribution  $\mu$.

\begin{center}
{\bf Acknowledgements}\end{center}

 The authors acknowledge the financial support of the DFG Forschergruppe ``Spectral analysis, asymptotic distributions and stochastic dynamics'' and SFB~701 ``Spectral structures and topological methods in mathematics''. Yu.K. was partially supported by the DFG research grant TR120/12-1. E.L. acknowledges the financial support of the SFB~611, Bonn University.

\end{document}